\documentclass{amsart}
\calclayout
\usepackage{graphicx} 
\usepackage{hyperref}
\usepackage{amssymb}
\usepackage{tikz-cd}
\newcommand{\ad}{\operatorname{ad}}
\newcommand{\SL}{\operatorname{SL}}
\newcommand{\C}{\mathbb{C}}
\newcommand{\Id}{\operatorname{Id}}
\newcommand{\tr}{\operatorname{tr}}
\newcommand{\Ch}{\mathrm{Ch}}
\newcommand{\Chs}{\underline{\mathrm{Ch}}}

\newcommand{\QC}{\mathcal{QC}}
\newcommand{\End}{\operatorname{End}}
\newcommand{\Rep}{\operatorname{Rep}}

\newtheorem{remark}{Remark}
\newcommand{\mf}{\mathfrak}
\newcommand{\ot}{\otimes}
\title{encyclopedia.tex}
\address{School of Mathematics, University of Edinburgh, Edinburgh, UK}
\email{D.Jordan@ed.ac.uk}
\author{David Jordan}

\date{May 2023}

\begin{document}
\title{Quantum character varieties}
\begin{abstract}
    In this survey article for the Encyclopedia of Mathematical Physics, 2nd Edition, I give an introduction to quantum character varieties and quantum character stacks, with an emphasis on the unification between four different approaches to their construction.
\end{abstract}

\maketitle

\section{Introduction}
The purpose of this article is to introduce the notion of a character variety, to explain its central role in mathematical physics -- specifically gauge theory  -- and to highlight four different approaches to constructing its deformation quantization, what are colloquially called ``quantum character varieties".  Rather remarkably, each distinct mechanism for quantization is motivated in turn by a distinct topological observation about the topology of surfaces and hence the geometry of classical character varieties.

The term ``character variety" commonly refers to a moduli space of $G$-local systems on some topological space $X$, equivalently of homomorphisms $\rho:\pi_1(X)\to G$.  In this article we discuss several different models for this moduli space -- the ordinary, framed, and decorated character varieties -- as well as their common stacky refinement, the character stack.  The distinction between these different models becomes important due to the presence of stabilizers and singularities in the naively defined moduli problem.

The term ``quantum character variety" typically refers to any of the following non-commutative deformations of character varieties in the case $X=\Sigma$ is a real surface:
\begin{enumerate}
	\item \textbf{Skein algebras} of surfaces yield deformation quantizations of the algebra of functions on ordinary character varieties of surfaces.
	\item \textbf{Moduli algebras} attached to ciliated ribbon graphs, yield deformation quantizations of the algebra of functions on framed character varieties of surfaces with at least one boundary component.
	\item \textbf{Quantum cluster algebras} associated to marked and punctured surfaces quantize a cluster algebra structure on the decorated character variety.
	\item \textbf{Factorization homology} of surfaces, with coefficients in the ribbon braided tensor category of representations of the quantum group yields linear abelian categories quantizing the category of quasi-coherent sheaves on  character stacks. 
\end{enumerate}

Each of the above constructions depends on a non-zero complex parameter $q$, and reproduces the corresponding ``classical" or ``undeformed" character variety upon setting $q=1$.  By definition the classical character variety depends on the underlying space $X$ only through its fundamental group -- in particular, only up to homotopy.  The quantum character varieties are more subtle invariants, which depend on the homeomorphism type of the manifold in a way detected only when $q\neq 1$.

In this survey we will recount the classical and quantum versions of each of these four constructions, outline their relation to one another, and how their study relates to super-symmetric quantum field theory and the mathematical framework of topological field theory.  The topological input for each construction above is a surface, however in each case there are natural extensions to 3-manifolds -- some more developed than others -- which we will also discuss.

\subsection{Flat connections}

One of the most fundamental notions in gauge theory is that of a principal $G$-bundle with connection.  A principal $G$-bundle $E$ over a space $X$ consists of a map $\pi:E\to X$, together with a free $G$-action on $E$ which preserves fibers and makes each fiber of $\pi$ into a $G$-torsor, so that $E/G=X$.  A connection on $E$ is a 1-form $A\in \Omega^1(X,\ad(E))$ valued in the adjoint bundle,
\[
\ad(E) = \mathfrak{g}\times_G E = (\mathfrak{g}\times E)/G.
\]

Among all connections on $E$ are distinguished the \emph{flat connections} $A$ which satisfy the flatness equation $dA + [A,A]=0$.  Equivalently $A$ is flat if the parallel transport along some curve $\gamma$ defined by $A$ depends on $\gamma$ only up to homotopy of paths.

\subsection{Local systems}
From a principal $G$-bundle with flat connection we may extract the more combinatorial data of a principal $G$-bundle $E$, together with parallel transport isomorphisms $\nabla_\gamma:E_x\to E_y$ along the homotopy class of paths connecting $x$ to $y$.  Such a pair $(E,\nabla)$ is called a $G$-local system.  The choice of a base point $x$ in $X$, and a trivialisation of $E$ at $x$ reduces the data of a $G$-local system to that of a group homomorphism $\pi_1(X)\to G$.  Changes of basepoint and changes of framing are both implemented by conjugation in $G$; hence two such homomorphisms represent the same $G$-local system if, and only if, they are related by post-composition with conjugation in $G$.

\subsection{Appearances in physics}
Classical character varieties appear very naturally in classical gauge theories such as Yang-Mills and Chern-Simons theory, in which a gauge field is by definition an adjoint valued 1-form, and the classical equations of motion involve the curvature of the connection -- in particular, in Chern--Simons theory in 3 dimensions, the critical locus is precisely the space of flat connections.

Quantum character varieties/stacks play a similar role in super-symmetric quantum field theories.  Some notable examples are:
\begin{enumerate}
\item Quantum character varieties -- specifically in their skein-theoretic incarnation -- describe topological operators, known as Wilson lines, in the quantization of Chern-Simons theory \cite{Witten-Jones}.  The parameter $\hbar=\log q$ appears as the quantization parameter.  
\item It is expected that state spaces of 3-manifolds, and categories of boundary conditions on surfaces, both attached to the Kapustin--Witten twist \cite{kapustin2006electric} of 4D $\mathcal{N}=4$ super Yang-Mills may be described via skein modules, and skein categories, respectively.  Here, the parameter $\Psi=\log q$ identifies with the twisting parameter in Kapustin--Witten's construction.
\item Coulomb branches of 4d $\mathcal{N}=2$ theories of class S, compactified on a circle, are naturally described by character varieties.  In this case, the deformation parameter $q$ appears as the exponentiated ratio of $\Omega$-deformation parameter for a pair of commuting circle actions coming from $SO(2)\times SO(2)\subset SO(4)$.  \cite{GMN, GMN2,hollands2016spectral,tachikawa2015skein}  In more mathematical terms, Gaiotto has proposed quantizations of Coulomb branches (hence of character varieties) via $\C^\times$-equivariant cohomology.  This has been discussed in \cite{schrader2019k}.
\end{enumerate}


\subsection{Acknowledgements}
In preparing this survey article I have benefited from the inputs of a number of people.  I particularly thank David Ben-Zvi, Adrien Brochier, Fran{\c c}ois Costantino, Charlie Frohman, Andy Neitzke, Thang Le, Gus Schrader, and Alexander Shapiro for fact-checking and helping to find complete references.  I also thank the Editor Catherine Meusberger, whose detailed comments and suggestions considerably improved the exposition.

\section{Classical character varieties}\label{sec:classical}
Each of the four quantization prescriptions highlighted in the introduction emerges from first understanding the relevant classical moduli space: its geometry, and its Poisson structure. When viewed from the correct perspective, each classical moduli space admits a very natural deformation quantization.  For this reason, although the article is focused on quantum character varieties, a significant portion is dedicated to reviewing the classical structures.

Let us begin by recalling in more detail the precise construction of the framed character variety, the ordinary character variety (the word ``ordinary" is non-standard, and appears here and throughout merely for emphasis), and the decorated character variety. We turn then to the ordinary and decorated character stacks, and then finally to enumerating the relations between them.  Along the way we will introduce the Poisson brackets which will be the quasi-classical limits of the quantum constructions.

\subsection{Framed character variety}  Given a group $G$ and a compact manifold $X$ with basepoint $x$, the framed character variety $\Ch_G^{fr}(X)$ is an affine algebraic variety which parameterises pairs $(E,\eta)$, where $E$ is a $G$-local system, and $\eta:E_x\to G$ is a trivialisation of the fiber $E_x\cong G$ (equivalently, this is the data of a single point $e\in E_x$ which becomes the identity element in $G$ under the framing).  In more concrete terms, we may identify the framed character variety with the set of representations $\pi_1(X)\to G$, where we do not however quotient by the $G$-action.  This alternative description makes it clear that $\Ch_G^{fr}(X)$ is indeed an algebraic variety: choosing any presentation of $\pi_1(X)$ with $m$ generators and $n$ relations identifies $\Ch_G^{fr}(X)$ with a closed subvariety of the affine variety $G^m$ defined by the $n$ relations, regarded as $G$-valued equations.  As there is no a priori reason for this closed subvariety to be equidimensional, the framed character variety will typically not be smooth.

An important special case is the framed character variety of a surface $\Sigma_{g,r}$ of genus $g$ with $r\geq 1$ punctures.  Since $\pi_1(\Sigma_{g,r})$ is the free group on $2g+r-1$ generators, we have that
\[
\Ch_G^{fr}(\Sigma_{g,r}) = G^{2g+r-1}.
\]

\subsection{Ordinary character variety}
The (ordinary) character variety $\Ch_G(X)$ is defined as the GIT quotient\footnote{with trivial stability condition, equivalently the ``categorical quotient" or ``coarse moduli space"} of the framed character variety by the $G$-action by conjugation.  By definition, this means that the character variety is an affine variety defined as
\[\Ch_G(X) = \operatorname{Spec}(\mathcal{O}(\Ch_G^{fr}(X))^G),\]
the spectrum of the sub-algebra of $G$-invariant functions on the framed character variety.

More geometrically, the character variety so defined parameterises \emph{closed $G$-orbits} in the framed character variety.  The map sending a point of the framed character variety to the closure of its $G$-orbit defines a surjection $\Ch_G^{fr}(X)\to \Ch_G(X)$.

\subsection{Decorated character variety}
For simplicity we will restrict now to the case where $X=\Sigma$ is a compact surface possibly with boundary.  An important enhancement of the notion of a $G$-local system is that of a decorated local system.  In a series of three highly influential papers \cite{FG06,FG09a,FG09b}, Fock and Goncharov established that the moduli space of decorated local systems, known as the decorated character variety, has an open locus admitting the geometric structure of a cluster variety, and they exploited this structure to define its quantization (discussed in Section \ref{sec:quantum}).

\begin{figure}[h]
    \includegraphics[height=1in]{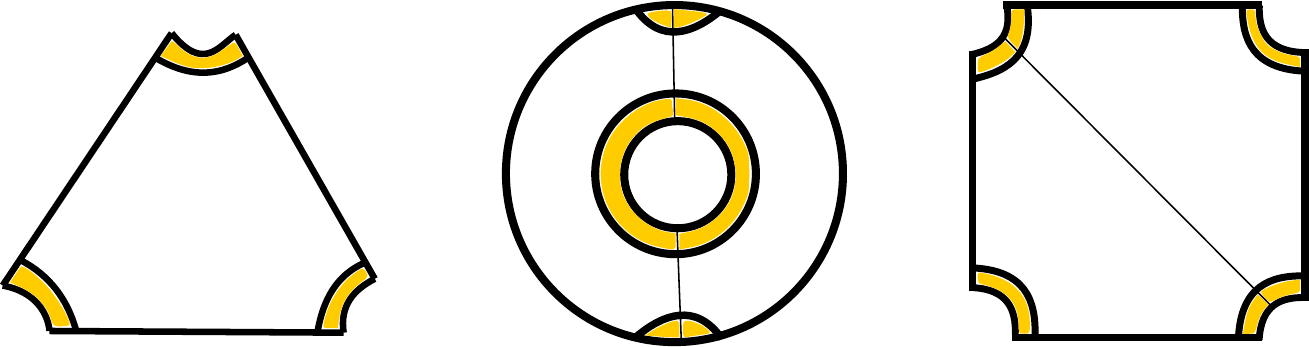}
    \caption{Three marked surfaces: the ``triangle", a disk with three contractible T-regions (indicated by shading), the ``punctured disk with two marked points", an annulus with one annular $T$-region and two contractible $T$ regions, and ``the punctured torus", the torus (opposite edges are identified) with a disk at the corner removed and an annular $T$-region around the resulting boundary circle.  The thin lines on the latter two depict a triangulation. }\label{fig:clusters}
\end{figure}

For the construction of decorated character varieties we fix in addition to the group $G$ a Borel subgroup $B$ and let $T=B/[B,B]$ denote the quotient of $B$ by its unipotent radical.  A $G$-$B$-$T$-coloring of a surface $\Sigma$ consists of a partition of $\Sigma$ into three-sets $\Sigma= \Sigma_G\sqcup \Sigma_B\sqcup \Sigma_T$, where $\Sigma_G$ and $\Sigma_T$ are open and $\Sigma_B=\partial \Sigma_G \cap \partial \Sigma_T$ is closed (see Figure \ref{fig:clusters} for some examples).  We will say that a ``marked point" on a decorated surface $\Sigma$ is a $T$-region which contracts onto an interval in the boundary of $\Sigma$, while a ``puncture" is a $T$-region contracting onto an entire boundary component of $\Sigma$.  We will call a connected decorated surface all of whose $T$ regions are of those two types a ``marked and punctured surface".  We note that a marked and punctured surface necessarily has a unique $G$-region.

\begin{remark}   
In Fock and Goncharov's original work, and most works which follow them, the marked points and punctures are indeed regarded as a finite set of points of $\Sigma$ lying in the boundary and interior of $\Sigma$, respectively, rather than as two-dimensional regions contracting to the boundary, as we have described above.  However when one unpacks the data they attach to punctures and marked points, one sees that it expresses quite naturally in the framework of decorated surfaces, and in particular the resulting notion of decorated local system, and hence the decorated character variety defined in either convention is identical.

For the topological field theory perspective it is important to ``zoom in" on these points and see them as one-dimensional ``defects" between adjoining 2-dimensional regions (the $G$- and $T$-regions discussed above).  For example, the ``amalgamation" prescription of Fock and Goncharov -- by which one glues together charts on each triangle of a triangulation to obtain a chart on the decorated character variety -- is just an instance of the excision axiom in factorization homology.
\end{remark}

The decorated character variety is a moduli space parameterising triples $(E_G,E_B,E_T)$, where $E_G$ and $E_T$ are $G$- and $T$-local systems over $\Sigma_G$ and $\Sigma_T$, respectively, and where $E_B$ is a reduction of structure over $\Sigma_T$ of the product $E_G\times E_T$ restricted there.  By a reduction of structure, we mean a $B$-sub-local system of the $B$-space $E_G\times E_T$.  At each point of the curve $\Sigma_B$ this is simply the specification of a $B$-orbit $B\cdot (g,t)\subset G\times T$, for some $(g,t)\in G\times T$ equivalently a point of $(G\times T)/B \cong G/N$; however the local system $E_G$ and $E_T$ can twist as we traverse loops in the surface, so that monodromy around punctures can introduce multiplication by some fixed element $t\in T$.  In other words, the monodromy around punctures need only preserve the underlying flag $\overline{\mathcal{F}}\in G/B$ obtained as the image of $\mathcal{F}\in G/N$ under the projection $G/N\to G/B$.

As with the ordinary character variety, one may introduce a framed variant of the character variety, by requiring additional data of a trivialisation of $E$ at a sufficiently rich system of basepoints.  It suffices to assume there is one basepoint in each connected region of $E_G$ and of $E_T$, and to require a trivialisation of the fiber there: in this case one can show that the resulting groupoid is rigid, so that the framed decorated character variety parameterising this data really is a variety.  We may then construct the GIT quotient of $\Ch_G^{fr,dec}(\Sigma)$ by $G^r\times T^s$, where there are $r$ basepoints in $\Sigma_G$ and $s$ basepoints in $\Sigma_T$.  It is an exercise analogous to that for the ordinary character variety to see that this construction is in fact independent of the chosen basepoints up to unique isomorphism.

In passing from the framed decorated character variety to the decorated character variety, we have quotiented by a typically non-free action of $G^r$, which has the consequence that the GIT quotient is typically singular.  Inspired by the Penner coordinate system \cite{penner1987decorated} on decorated Teichm\"uller space, Fock and Goncharov found a remarkable open subvariety of the framed and decorated character variety, on which the $G^r$-action is actually free, hence the GIT quotient there is smooth.  Moreover, they produced a remarkable set of ``cluster charts" in the framework of cluster algebras as had been introduced in \cite{ClustersI,ClustersII, ClustersIII,ClustersIV, gekhtman2002cluster}.  We recall the rudiments of their construction here.

Fock and Goncharov defined a family of subsets $U_\alpha$ on the framed decorated character variety of a marked and punctured surface with three remarkable properties:
\begin{enumerate}
    \item The $G$-action on $\Ch_G^{dec,fr}(\Sigma)$ restricts to a \emph{free} $G$-action on $U_\alpha$, and
    \item The quotient $U_\alpha/G$ is a torus $(\C^\times)^k$, for some explicitly given $k$.
    \item The transition maps between charts $U_\alpha/G$ and $U_\beta/G$ take the form of a ``cluster mutation", an explicitly given birational transformation between tori.
\end{enumerate}
The decorated character variety is defined as the union over the charts $\alpha$ of the $U_\alpha/G$; it is a subvariety of the decorated character stack, which although a union of affine charts, is typically not affine.

We emphasise moreover that the remaining $T^s$-action on each $U_\alpha/G$ is still not free; different ways of treating this non-free action, as well as specifying the monodromy of the $G/N$-fibers around punctures, leads to various related formulations of decorated character varieties.  We will not recall their complete definitions in this survey article.

The collection of opens $U_\alpha/G$ form what is called a cluster structure: briefly, this means that one has a combinatorial prescription to reconstruct the union of the charts $U_\alpha/G$ by starting from a single initial chart $U_0/G$ -- together with its coordinates and its Poisson structure $U_0/G$ is called a seed -- and successively adding in new charts, glued via the cluster mutation.  To construct the seed, one may choose a triangulation of the surface $\Sigma$, which must be compatible with the decoration, in that each vertex of the triangulation should be some framed basepoint lying in $\Sigma_T\cap \partial \Sigma$.  Each triangle contributes a torus, and one combines the different tori together along the edges of the triangulation via a process called ``amalgamation".  The coordinates and the Poisson structure are also specified in this construction.   The end result is also an algebraic torus, whose coordinates are indexed by the vertices of a quiver $\Gamma$, and whose Poisson bracket is by construction log-canonical; the rest of the cluster charts $U_\alpha/G$ are indexed by mutated graphs, and the precise form of the cluster mutation $U_\alpha/G\to U_\beta/G$ is encoded in this graph.

We will see in Section \ref{sec:quantum} that the explicit combinatorial description of cluster charts underpins an equally explicit Poisson structure and its canonical quantization to a quantum cluster algebra.

\subsection{The character stack}

In each of the above approaches, the presence of stabilisers inhibits a naive definition and introduces complications: for character varieties, we retreat to a moduli space of closed $G$-orbits, and for decorated character varieties, we retreat to an open locus where the $G$-action becomes free.  On the other hand, in both the ordinary and decorated case, our presentation has passed implicitly through a more universal notion of a character stack.

Without recounting completely the framework of stacks, we will recall only that certain moduli problems -- including both that of classifying ordinary and decorated local systems up to isomorphism -- admit the structure of an Artin stack: this simply means that they may be presented as the group quotient of an algebraic variety by a reductive algebraic group.  In fact, such a presentation is typically only required locally, but in our case we have a global such description.

One studies stacks algebraically via their locally presentable abelian categories of quasi-coherent sheaves, in particular we may consider the locally presentable abelian categories,
\[
\QC(\Chs_G(X)), \quad \QC(\Chs_G^{dec}(X)).
\]

As for any stack, $\QC(\Chs_G(X))$ carries a distinguished object called the \emph{structure sheaf} $\mathcal{O}$.  It follows from basic definitions that the algebra of functions on the ordinary character variety is isomorphic to $\End(\mathcal{O})$.  On the other hand, we have a pullback square,
\[
\begin{tikzcd}
\Ch_G^{fr}(X) \arrow{r}{} \arrow[swap]{d}{} & \Chs_G(X) \arrow{d}{g} \\
pt \arrow{r} & pt/G.
\end{tikzcd}
\]
Hence we may describe structures on $\Chs_G(X)$ as $G$-equivariant structures on $\Ch_G^{fr}(X)$, and conversely we recover $\Ch_G^{fr}(X)$ by forgetting this equivariance.

The relation between decorated character stacks and decorated character varieties of Fock and Goncharov is somewhat more complicated.  Each cluster chart $U_\alpha$ of the cluster variety $\Ch_G^{dec}(\Sigma)$ defines an object $\mathcal{O}_\alpha \in \QC(\Ch_G^{dec}(\Sigma))$: this is just the sheaf of functions which are regular on $U_\alpha$.  We have that $\End(\mathcal{O}_\alpha)$ is a ring of Laurent polynomials (i.e. functions on the corresponding torus $U_\alpha$), that the full subcategory generated by $\mathcal{O}_\alpha$ is indeed affine, and finally that the transition maps between the different $U_\alpha$'s define exact functors between these subcategories, which are written explicitly as cluster transformations.

A crucial feature of classical character stacks is that they fit into the framework of fully extended TQFT.  For this we briefly recall that there is a ``topological operad" $E_n$ which encodes the embeddings of finite disjoint unions of disks $\mathbb{R}^n$ inside one another.  An $E_n$-algebra is an algebraic structure governed by $E_n$.  $E_n$-algebras may be regarded as the ``locally constant" or ``topological" specialisation of the notion of a factorisation algebra, which in physical terms captures the structure of local observables in a quantum field theory, and the condition of being locally constant is a consequence when the QFT is topological.  Most relevant to our discussion are the examples that $E_1$-, $E_2$- and $E_3$-algebras in the bi-category of linear categories are monoidal, braided monoidal, and symmetric monoidal categories, respectively.  We refer to \cite{ayala2015factorization} for more technical details, and to \cite{brochiernotes} for a gentle exposition.

One may regard the symmetric monoidal category $\Rep(G)$ of representations of the reductive algebraic group $G$ as an $E_n$-algebra, for any $n$, in the symmetric monoidal 2-category of categories.  The collection of such $E_n$-algebras form an $n+1$-category, and the factorisation homology defines a full extended $n$-dimensional topological field theory valued in the symmetric monoidal 2-category of categories.  We have an equivalence of categories \cite{BZFN}:
\begin{equation}\label{eqn:FHom}
\QC(\Chs_G(X)) \simeq \int_X \Rep(G),
\end{equation}
where the integral notation on the righthand side denotes the factorization homology functor defined in \cite{ayala2015factorization}.  This equivalence and its consequences are discussed in greater detail in Section \ref{sec:quantum} below. 
 A similar TFT description can be given for decorated character stacks using stratified factorization \cite{AFT}.  For now, we note that it is not possible for either the ordinary or decorated character varieties to fit into the fully extended framework: if a manifold is given to us by gluing simpler manifolds, we can indeed build any $G$- or $G,B,T$-local system on it by gluing local systems on each piece, however there will be automorphisms of the glued local system which are not the disjoint product of automorphisms on each piece.  This simple observation prevents the ordinary and decorated character varieties from satisfying the gluing compatibilities satisfied by their stacky enhancements.

\section{Poisson brackets}\label{sec:Poisson}

Each of the framed, ordinary, and decorated character varieties, as well as the ordinary and decorated character stacks carry canonically defined Poisson brackets, which form the quasi-classical -- or leading order -- degenerations of the quantizations constructed in the next section.  We review these here.

\subsection{The Atiyah-Bott symplectic structure}

For reductive group $G$, Atiyah and Bott constructed in \cite{atiyah1983yang} a symplectic form on the smooth locus of the moduli space of flat $G$-connections on $\Sigma$.  Given a $G$-bundle $E$ with a flat connection $A$ on $\Sigma$, the tangent space to $E$ consists of sections of the associated adjoint bundle $\ad(E) = E\times_G \mathfrak{g}$ over $\Sigma$.  We denote by $\kappa$ the Killing form on $\mathfrak{g}$. Given a pair $\chi_1, \chi_2$ of sections of $\ad(E)$, regarded as tangent vectors to $E$, we define their symplectic pairing by the formula:

\[
\omega(\chi_1,\chi_2) = \int_{\Sigma} \kappa(\chi_1\wedge \chi_2).
\]

We remark that this symplectic form arises naturally in the Chern-Simons action term for a 3-manifold of the form $\Sigma\times I$.  The skew-symmetry, non-degeneracy and closedness of $\omega$ all follow relatively easily from the definition.

\subsection{The Goldman bracket on the ordinary character variety}
The character variety carries a combinatorial analog of the Atiyah--Bott symplectic structure, due to Goldman \cite{goldman1984symplectic} and Karshon \cite{karshon1992algebraic}, which can be defined purely algebraically without appeal to analysis.  The fundamental group $\pi_1(\Sigma)$ of a closed surface has as its second group cohomology $H^2(\pi_1(G),\mathbb{C}) = \mathbb{C}$.  The tangent space to a given representation $\rho:\pi_1(\Sigma)\to G$ identifies with $H^1(\pi_1(\Sigma),\mathfrak{g}),$ where $\pi_1(\Sigma)$ acts on $\mathfrak{g}$ through the conjugation action.  In analogy with the Atiyah-Bott symplectic structure we obtain a skew pairing,
\[H^1(\pi_1(\Sigma),\mathfrak{g})\times H^1(\pi_1(\Sigma),\mathfrak{g}) \xrightarrow{\kappa\circ \cup} H^2(\pi_1(\Sigma),\mathbb{C}) = \mathbb{C},\]
which again defines a symplectic form on the character variety.  The corresponding Poisson bracket is known as the Goldman bracket.

\subsection{The Fock-Rosly Poisson structure on the framed character variety}
While both the Atiyah-Bott and Goldman symplectic forms are clearly very natural and general, they don't lead immediately to explicit formulas for the Poisson brackets of functions on the character variety.  In the case of the framed character variety of surfaces with at least one boundary component, a much more explicit reformulation was given by Fock and Rosly in \cite{fock1998poisson}.

First we recall that $\pi_1(\Sigma_{g,r})$ is the free group on $2g+r-1$ generators, where $\Sigma_{g,r}$ denotes a surface of genus $g$ with $r$-punctures.  Hence the framed character variety is simply the product $G^{2g+r-1}$.  The Poisson bracket between between functions $f$ and $g$ of a single $G$ factor is given by the Poisson bivector,
\begin{equation}\label{eqn:STS}
\pi_{STS}=\rho^{ad,ad}+t^{r,l}-t^{l,r}.
\end{equation}
Here we denote with superscripts $r,l, ad$ the right, left and adjoint vector fields on $G$ determined by a Lie algebra element, we let $r\in (\mf g)^{\ot 2}$ denote the classical $r$-matrix, and $\rho$ and $t$ its anti-symmetric and symmetric parts.  The bivector $\pi_{STS}$ induces on $G$ a Poisson structure which has been introduced by Semenov--Tian--Shansky~\cite{STS94}.

The Poisson bracket between functions $f$ and $g$ of the $i$th and $j$th factor is given by $\pi_{ij}$, where
\[
\pi_{ij}=
\begin{cases}
\pm (r^{ad,ad}) & \text{if $i,j$ are $\pm$ unlinked}\\
\pm (r^{ad,ad}+2t^{r,l}) & \text{if $i,j$ are $\pm$ linked}\\
\pm (r^{ad,ad}-2t^{r,r}+2t^{r,l}) & \text{if $i,j$ are $\pm$ nested }
\end{cases}
\]

In total we have,
\begin{equation}
\pi = \sum_i \pi_{STS}^{(i)}+\sum_{i<j} (\pi_{ij}-\pi_{ji}),
\end{equation}
where now $\pi_{ij}$ is a 2-tensor acting on the $i$th component of the first factor and the $j$th component of the second factor, given by the formula above.

The appearance of the classical $r$-matrix in the Fock-Rosly Poisson bracket foreshadows the role of quantum groups and quantum $R$-matrices in the deformation quantization of Section \ref{sec:quantum}.

\subsection{The Fock-Goncharov cluster Poisson structure on the decorated character variety}

Recall that the decorated character variety of a marked and punctured surface contains a system of open charts, each isomorphic to a torus $(\mathbb{C}^\times)^k$ for some $k$ depending on $G$ and on the decorated surface.  Each chart $U_\alpha$ carries a ``log canonical" Poisson bracket:
\[
\{x_i,x_j\} = a_{ij}x_ix_j,
\]
where $A = (a_{ij})$ is the adjacency matrix of the quiver attached to $U_\alpha$. It is called log-canonical because the formal logarithm of the generators $x_i$ satisfy $[\log(x_i),\log(x_j)] = a_{ij}$, which resemble Heisenberg's ``canonical commutation relations".  The birational transformations given by cluster mutation intertwine the different Poisson brackets on each chart, so that they glue together to a globally defined Poisson bracket on the cluster Poisson variety.

\subsection{Shifted symplectic structures and character stacks}

The moduli problem given by the character stack can be phrased in terms of classifying spaces, and this allows a universal construction of the Atiyah--Bott/Goldman symplectic structure on character varieties, due to Pantev, To\"en, Vaquie, and Vezzosi \cite{pantev2013shifted}, see \cite{CALAQUE} for an exposition.  Specifically, there exists a classifying stack $BG$, such that for any surface $\Sigma$, we have:
\[
\Chs_G(\Sigma) = \operatorname{Maps}(\Sigma,BG),
\]
where $\operatorname{Maps}$ here denotes the stack of locally constant maps from a topological space into an algebraic stack, for instance as obtained by regarding $\Sigma$ as presented by a simplicial set. 

The classifying space $BG$ has as its tangent and cotangent complexes the Lie algebra $\mathfrak{g}$ and its dual $\mathfrak{g}^*$ in homological degrees 1 and -1, respectively.  Hence the Killing form gives an isomorphism from the tangent bundle to the 2-shifted cotangent bundle: such a structure (satisfying some additional properties which follow from those of the Killing form) is known as a 2-shifted symplectic structure.  One may then transgress the 2-shifted symplectic structure on $BG$ through the mapping stack construction to give an ordinary -- or 0-shifted -- structure on the character stack.    Remarkably, the descent of this symplectic structure to the smooth part of the character variety recovers the Atiyah-Bott/Goldman symplectic structure.  Hence the PTVV structure on the character stack may be regarded as a stacky version of the Atiyah--Bott/Goldman construction.

\subsection{Hamiltonian reduction interpretation}
The framework of Hamiltonian reduction gives a natural and very general procedure to pass from the framed character variety of the surface $\Sigma_{g,1}$ to the ordinary character variety of the closed surface $\Sigma_g$.

Let $\mu: \Ch_G^{fr}(\Sigma_{g,1})\to G$ denote the map which sends the $G$-local system to its monodromy around the unique boundary component (we assume for simplicity that the basepoint is on the boundary to ensure this map is well-defined, not only up to conjugation).  Sealing the boundary component means imposing $\mu=\Id$, and forgetting the framing means quotienting by the $G$-action.  Hence we have,
\[
\Ch_G(\Sigma_g) = \mu^{-1}(\Id) / G.
\]
Formulas such as the above are common in the theory of \emph{Hamiltonian reduction}, where one interprets $\mu$ as a \emph{moment map} for a Hamiltonian action of the group $G$ on some phase space -- in this case the phase space is $\Ch^{fr}_G(\Sigma_{g,1})$. The target of the moment map is typically $\mathfrak{g}^*$ rather than $G$, however Alekseev, Kosmann-Schwarzbach, Malkin, and Meinrenken developed in \cite{alekseev1998lie, alekseev2002quasi}, following \cite{lu1991momentum}, a formalism of ``quasi-Hamiltonian" $G$-spaces, which feature ``group-valued moment maps" valued in $G$ rather than $\mathfrak{g}^*$.  The Hamiltonian reduction of a quasi-Hamiltonian $G$-space is symplectic, and recovers to the Atiyah-Bott/Goldman symplectic structure on the closed surface.

\section{Quantum character varieties}\label{sec:quantum}
In this final section, let us recount the four most well-known constructions of quantum character varieties.

\subsection{The moduli algebra quantization of the framed character variety}
The Fock-Rosly Poisson structure on the framed character variety of $\Sigma_{g,r}$, with $r\geq 1$, admits a highly explicit and algebraic quantization, introduced independently by Alekseev, Grosse, and Schomerus \cite{AGSI,AGSII} and by Buffenoir and Roche \cite{buffenoir1995two}. This  quantization was orginally called the ``moduli algebra", and is often referred to subsequently as the ``AGS algebra".  The constuction also goes under the names ``combinatorial Chern-Simons theory" and ``lattice gauge theory", and has been studied in many different contexts, see e.g. \cite{roche2002trace}, \cite{MeusbergerSchroers}, \cite{meusburger2017kitaev}.
The starting point is that the classical $r$-matrices appearing in the Fock-Rosly Poisson bracket have well-known quantizations into quantum $R$-matrices, which themselves describe the braiding of representations of the quantum group.  Artful insertion of quantum $R$-matrices in place of the classical $r$-matrices gives the deformation quantization of framed character varieties: let us now describe their construction in more detail.

In the case $g=1,r=1$, the Fock-Rosly Poisson bracket is identical to the Semenov-Tian-Shansky Poisson bracket $\pi_{STS}$ from Equation \eqref{eqn:STS}. 
 Replacing classical $r$-matrices with quantum $R$-matrices leads to a deformation quantization $\mathcal{F}_q(G)$ of the algebra of functions\footnote{A point of disambiguation:  the reflection equation algebra $\mathcal{F}_q(G)$ here does not coincide with the so-called FRT algebra quantization, often denoted $\mathcal{O}_q(G)$, which quantizes instead the Sklyanin Poisson bracket on $G$.} on the group $G$, known as the \emph{reflection equation algebra}.  The name refers to the fact that for $G=\operatorname{GL}_N$, the commutation relations in this algebra are given\footnote{See, e.g, \cite{JordanWhite} for details about this notation.} via the ``reflection equation",
\begin{equation}\label{eqn:REA}
R_{21}A_1 R_{12} A_2 = A_2 R_{21}A_1R_{12},
\end{equation}
which appears as a defining relation in Coxeter groups of type $B$, and in mathematical physical models for scattering matrices in 1+1-dimension in the presence of a reflecting wall. 

A general surface with boundary may be presented combinatorially as a ``ciliated ribbon graph" -- essentially a gluing of the surface from disks.  According to this prescription, each edge contributes a factor of $\mathcal{F}_q(G)$ to the moduli algebra (which may be regarded as the quantum monodromy of a connection along that edge), and the cross relations between different edge factors are given by explicit formulas resembling the reflection equation, but with an asymmetries related to the unlinked, linked or nested crossings.

\subsection{The skein theory quantization of the ordinary character variety}
The skein algebra quantization is prefaced on an elegant graphical formulation of the functions on the classical $\SL_2$ character variety in terms of (multi-) curves drawn on the surface.  This is then deformed to a similar graphical calculus for curves drawn instead in the surfaces times an interval.  Skein algebras were independently introduced by Przytycki \cite{Przytycki} and Turaev \cite{Turaev}.  Following them, the vast majority of skein theory literature concerns the so-called Kauffmann bracket skein relations, a particular normalisation of the skein relations deforming the $\SL_2$-character variety.  In that same tradition, we will recall this special case first in detail, only outlining the extension to general groups, and indeed to general ribbon braided tensor categories afterwards.

Given a loop $\gamma: S^1\to X$, a $G$-local system $E$ over $X$, and a finite-dimensional representation $V$ of $G$, we have a polynomial function, $\tr_{\gamma,V}$, sending $E$ to the trace of the parallel transport along $\gamma$ of the associated vector bundle with connection $E\times_G V$.  The word ``polynomial" here means, for instance that $\tr_{\gamma,V}$ defines a $G$-invariant function on the framed character variety, hence a polynomial function on the GIT quotient. The skein theoretic approach to quantizing character varieties begins by describing this commutative algebra structure ``graphically", i.e. via the image of the curve $\gamma$ sitting in $\Sigma$, and then introducing a deformation parameter in the graphical presentation.

\begin{figure}
    \includegraphics[height=1in]{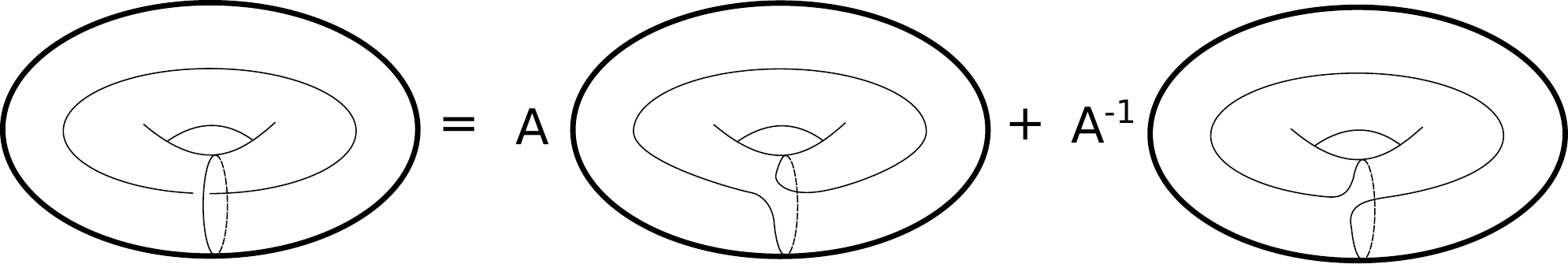}
\caption{The product $\tr_{\gamma_1,V} \cdot \tr_{\gamma_2,V}$ of two curves is given by their stacking in the $I$ direction on $\Sigma\times I$.  The skein relations express this as a linear combination of two new curves as depicted at the right.}\label{fig:skeinrelntorus}
\end{figure}

Let us consider a 3-manifold $M=\Sigma\times I$, the cylinder over some surface $\Sigma$.  We may depict the function $\tr_{\gamma,V}$ by drawing $\gamma$ with the label $V$ (one often projects onto the $\Sigma$ coordinate, to draw $\gamma$ as a curve with normal crossings on $\Sigma$).  We depict the product $\tr_{\gamma_1,V_1}\cdot \tr_{\gamma_2,V_2}$ of two such functions by super-imposing the drawings, as in Figure \ref{fig:skeinrelntorus}.  The resulting diagrams will of course develop additional crossings when multiplied, however a basic observation -- which is somewhat special to the case of $\SL_2$ -- is that one can use local ``skein relations" to resolve crossings.

\begin{figure}[h]
\includegraphics[height=1in]{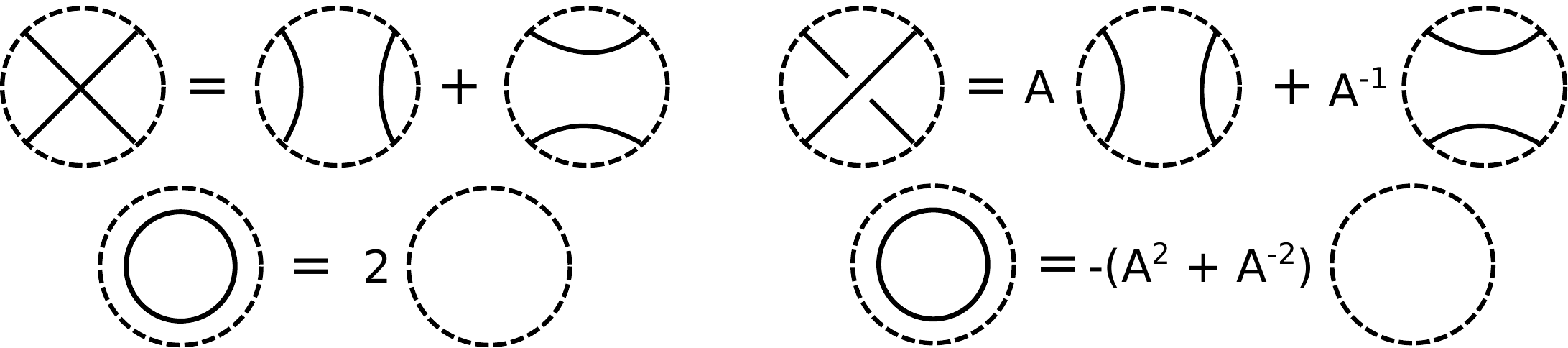}
\caption{At left, the equations $\tr(XY)+\tr(X^{-1}Y) = \tr(X)\tr(Y)$ and $\tr(\Id)=2)$ express graphically as skein relations; the top relation holds between any three curves which are identical outside of the dotted region, and differ as indicated within it, while the bottom relation states that any isolated unknot can be erased at the price of a factor of 2.  At right is stated the quantum deformation which depends on a parameter $A\in \mathbb{C}$.}\label{fig:skeinreln}
\end{figure}

To better understand this, consider the Cayley-Hamilton identity for a matrix $X\in \SL_2$:
\[ X + X^{-1} = \tr(X)\cdot\Id_2,\]
Multiplying by an arbitrary second matrix $Y$ and taking traces gives an identity
\begin{equation}\tr(XY)+\tr(X^{-1}Y) = \tr(X)\tr(Y).\label{eq:classical-skein}\end{equation} 
Suppose now that we have two paths $\gamma_1$ and $\gamma_2$ intersecting at a point $p \in \Sigma$, and that we consider the product $\tr_{\gamma_1,V_1}\cdot \tr_{\gamma_2,V_2}$, where $V_1=V_2=\mathbb{C}^2$ is the defining representation of $\SL_2$.  Then the identity \eqref{eq:classical-skein} implies a graphical simplification as depicted in Figures \ref{fig:skeinreln} and \ref{fig:skeinrelntorus}.

A fundamental observation of Przytycki and Turaev (independently) is that the graphical relation in Figure \ref{fig:skeinreln} can be naturally deformed by introducing coefficients into the relations which mimic the defining relations of the Jones polynomial, more precisely its normalisation known as the Kauffmann bracket.  Around the same time, Edward Witten had understood these same ``skein relations" to appear as fundamental relations between Wilson loops in a quantization of the Chern-Simons theory in 3-dimensions.

The $\SL_2$-skein algebra of a surface is therefore defined as the quotient of the vector space spanned by links embedded into $\Sigma\times I$, modulo local skein relations: two skeins are declared equivalent if they agree outside of some cylindrical ball $\mathbb{D}^2\times I$ and differ inside as indicated in the Figure.  An example of a skein relation is depicted in Figure \ref{fig:skeinrelntorus}: we imagine a thickening of $T^2$ to $T^2 \times I$, and a small cylindrical ball $D^2\times I$ around the intersection point indicated on the left hand side of the equation. Within the cylindrical ball, we replace the skein with the expressions in the right hand side of the skein relation in Figure \ref{fig:skeinreln}.  We impose such relations for all cylindrical balls throughout $\Sigma\times I$.

Turaev showed that the resulting relations are flat in the parameter $A$ -- more precisely he showed that the skein module\footnote{The name ``module" refers simply to the fact that the base ring may be taken to be $\mathbb{Z}[A,A^{-1}]$ (to allow specialisation) rather than a field.  However, when a 3-manifold has boundary, it skein module indeed becomes a module for the action by inserting skeins at the boundary.} is free as a $\mathbb{Z}[A,A^{-1}]$-module, and that the specialisation at $A=-1$ recovers the vector space of polynomial functions on the character variety of the surface\footnote{The careful reader may note a perhaps unexpected sign here, and in the skein relations at the right side of Figure \ref{fig:skeinreln}; this is a standard convention taken to ensure the Kauffman bracket skein module of a 3-manifold does not depends only on a choice of orientation. Although in the case of surfaces, the signs can be absorbed into the normalisation, we include the standard normalisation here for the sake of consistency. The parameter $A$ discussed here is a fixed square root of the parameter $q$ discussed elsewhere in the article.}.

Moreover, the skein module of $\Sigma\times I$ obtains an algebra structure by stacking skeins in the $I$-direction.  With respect to this stacking operation, the skein module becomes a non-commutative algebra whose $q=1$ specialisation is the algebra of functions on the classical ordinary character variety.  Due to the flatness in $q$, we obtain a Poisson bracket on the character variety, by setting:
\begin{equation}
\{f,g\} = \frac{f\cdot g - g\cdot f}{q-q^{-1}} \mod q-1.\label{eqn:degen}
\end{equation}
It was shown in \cite{bullock1999understanding} that this Poisson bracket agrees with the Atiyah--Bott/Goldman bracket and so the skein algebra is a deformation quantization of the character variety with its Atiyah--Bott/Goldman Poisson bracket.

The definition above extends naturally to define the $\SL_2$-skein module of an oriented 3-manifold: one considers the formal linear span of links in $M^3$, modulo the skein relations imposed in each cylindrical ball $D^2\times I$ embedded in the 3-manifold.  The $A=-1$ specialisation still coincides with the functions on the character variety of $M$ \cite{bullock1997rings}, however the skein module in general no longer carries an algebra structure, as there is no distinguished direction along which to stack the skeins.

The above discussion has been formulated for simplicity in the case of $\SL_2$-skeins, however it generalises naturally to any simple gauge group $G$, indeed to any ribbon tensor category.  The careful reader will have noticed that it sufficed in the $\SL_2$-case to consider only the defining two dimensional representation, and that we could reduce always to crossing-less diagrams.  For a general group $G$, or more generally for an arbitrary ribbon tensor category, this is no longer possible.  Skeins for a general group $G$ consist of embedded oriented ribbon graphs in the 3-manifold, together with a network of labels of representations of the quantum group (more genenerally objects of the ribbon braided tensor category) along each edge, and a morphism at every vertex, from the tensor product of incoming to outgoing edges.  We refer to \cite{cooke2019excision} or \cite{GJS} for details about the general construction, or to \cite{kuperberg1996spiders} or \cite{sikora2005skein} for early examples beyond the Kauffman bracket skein module.

An important ingredient in the definition is the Reshetikhin--Turaev evaluation map, which maps any skein on the cylindrical ball $\mathbb{D}^2\times I$ to a morphism from the tensor product of incoming labels along $\mathbb{D}^2\times \{0\}$ to the tensor product of the outgoing labels along $\mathbb{D}^2\times \{1\}$.   The skein module is defined as the span of all skeins, modulo the kernel of the Reshtikhin--Turaev evaluation maps, ranging over all embedded disks in $M$.

\subsection{The Fock--Goncharov quantum cluster algebra structure on the decorated character variety}
Recall that the cluster algebra structure on the decorated character varieties consists of a collection of open subsets, each identified via a coorinate system with a torus $(\C^\times)^r$, carrying a ``log-canonical" Poisson bracket preserved by the birational transformations.  There is a canonical quantization of such a torus given by introducing invertible operators satisfying $X_iX_j = q^{a_{ij}} X_jX_i$.  The Poisson bracket obtained from these relations as in equation \eqref{eqn:degen} recovers the log-canonical Poisson bracket.

An elementary and fundamental observation of Fock and Goncharov is that conjugation by the ``quantum dilogarithm" power series induces a birational equivalence between different quantum charts. This birational isomorphism lifts the cluster mutation taking place at the classical level to an birational isomorphism of the associated quantum tori.  Essentially by definition, the cluster quantization of the decorated character variety is the resulting collection of quantum tori, equipped with a preferred system of generators, known as cluster variables, related by quantum cluster mutations.  

In summary, Fock and Goncharov construct the quantization of the decorated character variety as a quantum cluster algebra.  The resulting algebraic, combinatorial, and analytic structures are of considerable interest and are heavily studied, however to survey them in proper depth would be beyond the scope of this survey article.  Instead we highlight several important papers following in this tradition: \cite{Le19, SS19, Tes,Ip18, GS19,FST08,schrader2017continuous},

\subsection{Quantum character stacks from factorization homology}

A basic ingredient in non-commutative algebra and gauge theory is the quantum group.  In categorical terms, the category $\Rep_q(G)$ is a ribbon braided tensor category, which $q$-deforms the classical category of representations of the algebraic group $G$.

Recall that the classical character stack is computed via factorization homology, as in Equation \eqref{eqn:FHom}.  The quantum character stack of a surface is defined, by fiat, by the equation,
\[
Z(\Sigma) = \int_{\Sigma} \Rep_q(G).
\]
Here, as in Section \ref{sec:classical} we are regarding $\Rep_q(G)$ with its braided monoidal structure equivalently as an $E_2$-category in the symmetric monoidal bi-category of categories, in order to define its factorization homology on surfaces.  We also note that the ribbon structure on $\Rep_q(G)$ equips it with an $SO(2)$-fixed point structure, so that it descends from an invariant of framed surfaces to one of oriented surfaces.

The construction via factorization homology opens up tools in topological field theory, extending the construction of quantum character stacks both \emph{down} to the level of the point, and \emph{up} to dimension 3.  The algebraic framework for the extended theory involves a 4-category denoted $\operatorname{BrTens}$, whose objects are braided tensor categories and whose higher morphisms encode notions of algebra, bimodules, functors and natural transformations, respectively, between braided tensor categories.  The braided tensor category $\Rep_q(G)$ defines a 3-dualizable object in $\operatorname{BrTens}$, and so according to the cobordism hypothesis it defines a fully extended framed 3-dimensional topological field theory.  The $SO(2)$-fixed structure descends this to an oriented topological field theory, and the resulting invariants of oriented surfaces coincide with the factorization homology functor $Z$ as defined above.

The construction by factorization homology has a modification where we allow a pair of braided tensor categories and a morphism between them labelling a codimension 1-defect, and from this data produce an invariant of bipartite surfaces.  Taking  $\Rep_q(G)$ and $\Rep_q(T)$, with $\Rep_q(B)$ labelling the defect, one obtains a quantum deformation of the decorated character stacks.

\subsection{Special structures at roots of unity}
Each flavour of quantum character variety discussed above exhibits special behavior when the quantum parameter $q$ is taken to be a root of unity:

\begin{enumerate}
    \item Skein algebras at root-of-unity parameter $q$ have large centers, over which the entire algebra is a finitely generated module, as first proved by Bonahon--Wong \cite{bonahon2016representations}.  Following them, numerous papers have studied the implications for the representation theory of the skein algebra, especially the determination of the Azumaya locus -- the points over the spectrum of the center for which the central reduction of the skein algebra is a full rank matrix algebra.  Most such results are only proved in the case of $\SL_2$-skeins, but are expected to hold more generally.  See \cite{frohman2019unicity,ganev2019quantum,karuo2022azumaya}
    \item AGS moduli algebras at a root of unity develop a large centre, and their Azumaya locus is known to contain the preimage $\mu^{-1}(G^\circ)$ of the big cell.  This follows easily from the Brown--Gordon/Kac--de Concini technique of Poisson orders \cite{brown2003poisson, de1991representations}:  the space $\mu^{-1}(G^\circ)$ is precisely the open symplectic leaf in the framed character variety, and BG/KdC method gives an isomorphism of the fiber over any two points in the same symplectic leaf.  See \cite{ganev2019quantum} for a precise formulation in the setting of AGS algebras.
    \item Quantum cluster algebras attached to decorated quantum character varieties exhibit parallel behavior at roots of unity:  chart by chart, the $\ell$th power of any cluster monomial is central, whenever $q^\ell=1$.  The cluster mutations respect this central structure, and lead to a quantum Frobenius embedding of the classical decorated character variety (of each respective type) into the corresponding quantum cluster algebra.  It is much more straightforward to determine the Azumaya locus on this setting, given the explicit control on the individual charts.  See \cite{FG09a}, and in greater generality \cite{nguyen2021root}.
    \item The Azumaya algebras described above are each an instance of \emph{invertibility}: indeed, an algebra $A$ is said to be invertible over its center of the algebra $A\otimes_{Z(A)} A^{op}$ is Morita equivalent to $Z(A)$.  
    A proof has been announced by Kinnear that quantum character stacks are invertible relative to classical character stacks in the strongest possible sense: the symmetric monoidal category $\QC(\Ch_G(\Sigma))$ acts monoidally on the quantum character stack, and the quantum character stack defines an invertible sheaf of categories for that action.  At the level of surfaces, this implies that the factorization homology category is an invertible sheaf of categories of quasi-coherent sheaves on the character stack.  More fundamentally, the assertion is that the fully extended 3D TQFT defined by quantum character stack construction at root-of-construction is invertible relative to the fully extended 4D TQFT determined by $\Rep(G)$.
\end{enumerate}

\subsection{Unifications of various approaches to quantum character varieties}

Recall that classically, each of the framed, ordinary, and decorated character stacks could be derived directly from the classical character stack under various geometric operations.  This applies equally well to the quantum character stacks.  We record the following unifications:

\begin{itemize}
    \item By a result of Cooke \cite{cooke2019excision}, the skein category is computed via the factorization homology of the category $\Rep_q(G)^{cp}$ of compact-projective objects.  The skein algebra is the algebra of endomorphisms of the empty skein object.
    \item In \cite{BBJ18a} Alekseev-Grosse-Schomerus algebras are recovered from the quantum character stack via monadic reconstruction techniques.  The quantum Hamiltonian reduction procedure is recast via monadic reconstruction in \cite{BBJ18b, safronov2019categorical}.
    \item The Fock--Goncharov quantum cluster algebra may be recovered via an open subcategory of the decorated quantum character stack: this is presently only proved in the $\SL_2$ case, in \cite{jordan2021quantum}.
    \item The Alekseev-Grosse-Schomerus algebras have also been recovered directly in skein-theoretic terms, via the so-called ``stated" \cite{Costantino-Le}, or ``internal" \cite{GJS} skein algebra construction.  See \cite{haioun2022relating}.
    \item The Bonahon--Wong approach to the representation theory of skein algebras involves describing skein algebras via quantum trace maps  (c.f. \cite{bonahon2011quantum}, \cite{kim2023sl2}, \cite{le2018triangular},).  This relationship has been studied in more physical literature under the term  non-abelianization in \cite{hollands2016spectral,GMN2, neitzke2020q,Neitzke_2020}.
\end{itemize}

\subsection{Quantum character varieties of 3-manifolds}
In contrast to the vast literature we have surveyed pertaining to quantum character varieties and character stacks, much less is currently known about the quantization of character varieties and character stacks of 3-manifolds.  Perhaps one reason for this, as discussed below, is that the very notion of quantization in the context of character stacks of 3-manifolds is different than in the case of surfaces: neither character stacks nor character varieties of 3-manifolds are symplectic/Poisson, but they are rather (-1)-shifted symplectic.  This implies a quantization theory of a different nature -- in particular as we discuss below, the skein module quantization of a 3-manifold is essentially never a flat deformation, this happens only when the character variety is a finite set of points.

Let us nevertheless highlight what is known and currently under investigation concerning each of the four perspectives on quantization.

\begin{enumerate}
    \item The AGS algebra attached to the surface $\Sigma_{g,1}$ acts naturally on the underlying vector space of the AGS algebra attached to $\Sigma_{0,g}$, which we should regard instead as attached to the handle body $H_g$ of genus $g$.  In \cite{GJS} it was established that the skein module of a closed oriented 3-manifold $M$ may be computed by choosing a Heegaard splitting,
    \[ M = H_g \cup_{\Sigma_g} H_g,\]
    where the splitting involves twisting by a choice of $\gamma$ in the mapping class group of $\Sigma_g$.
    \item This presentation of the skein module was used in \cite{GJS} to establish the finite-dimensionality of skein modules of closed 3-manifolds, as had been previously conjectured by Witten. Several recent papers \cite{carregaGeneratorsSkeinSpace2017,detcherryBasisKauffmanSkein2021, detcherryInfiniteFamiliesHyperbolic2021,gilmerKauffmanBracketSkein2007,gilmerKauffmanBracketSkein2018, detcherry2023kauffman,GJVY} are devoted to determining these dimensions in special cases, and Gunningham and Safronov have announced an identification of the skein module with the space of sections of a certain perverse sheaf introduced by Joyce as a quantization of the $(-1)$-shifted structure on the character variety regarded as a critical locus.
    \item Decorated character varieties of 3-manifolds (perhaps not by this name) were studied in the papers \cite{Dim,DGG}: one fixes an ideal triangulation of a hyperbolic 3-manifold, and directly computes a deformation quantization of the $A$-polynomial ideal using quantum cluster algebra techniques tuned to each ideal tetrahedron, regarded as a filling of a decorated $\Sigma_{0,4}$.
    \item It was established by Przytycki and Sikora \cite{Przytycki1997OnSA} that the $q=1$ specialisation of the skein module of $M$ indeed recovers the algebra of functions on the character variety of $M$. Echoing the Azumaya/invertibility at the level of surfaces, it is now expected that skein modules at root-of-unity parameters arise as the global sections of a line bundle on the classical character variety.  Two approaches to this result have been announced, one by Kinnear 
    using higher categorical techniques in parallel to the above discussion, and another by Frohman, Kania-Bartoszynska and Le, appealing to the Azumaya property for surfaces and invoking the structure of a (3,2) TQFT.
\end{enumerate}

\section{Further reading}
The following references may be helpful to a reader hoping to learn this subject in more detail:
\begin{itemize}
    \item ``Lectures on gauge theory and Integrable systems", \cite{Audin1997}.
    \item ``Quantum geometry of moduli spaces of local systems and representation theory" \cite{GS19}.
    \item \textit{Cluster algebras and Poisson geometry}, \cite{gekhtman2010cluster}.
    \item ``Factorization homology of braided tensor categories" \cite{brochiernotes}.
    \item GEAR Lectures on quantum hyperbolic geometry \cite{Frohmannotes}.
\end{itemize}

\bibliographystyle{plain}
\bibliography{bib}{}

\end{document}